\documentclass[times, 10pt,twocolumn]{article} 
\usepackage{latex8}
\usepackage{graphicx}
\usepackage{latexsym}

\title{The Fundamental Theorems of Interval Analysis
}

\author{M.H. van Emden and B. Moa \\
           Department of Computer Science\\
           University of Victoria, Canada
          }

\begin{document}
\vspace{-15mm}
\maketitle
\pagestyle{empty}
\vspace{-5mm}
\begin{abstract}
Expressions are not functions.
Confusing the two concepts or failing to define
the function that is computed by an expression
weakens the rigour of interval arithmetic.
We give such a definition and continue
with the required re-statements
and proofs of the fundamental theorems of interval arithmetic
and interval analysis.
\end{abstract}

\vspace{-5mm}
\newcommand{\cpp}{\hbox{{\tt C++}}}

\newtheorem{theorem}{Theorem}{}
\newtheorem{definition}{Definition}{}
\newtheorem{example}{Example}{}
\newtheorem{lemma}{Lemma}{}
\newtheorem{corollary}{Corollary}{}

\newcommand{\Emin}{\ensuremath{E_{\mbox{min}}}} 
\newcommand{\Emax}{\ensuremath{E_{\mbox{max}}}} 

\newcommand{\Rat}{\ensuremath{\mathcal{Q}}}
\newcommand{\Rea}{\ensuremath{\mathcal{R}}}
\newcommand{\ExtRe}{\ensuremath{\mathcal{R}^{++}}}
\newcommand{\Flpt}{\ensuremath{\mathcal{F}}}
\newcommand{\Intr}{\ensuremath{\mathcal{I}}}

\newcommand{\Alp}{\ensuremath{\mbox{\textbf{A}}}} 
\newcommand{\AlpStar} {\ensuremath{\Alp ^\ast}}

\newcommand{\A}{\ensuremath{\mathcal{A}}} 
\newcommand{\Bool}{\ensuremath{\mathcal{B}}} 
\newcommand{\C}{\ensuremath{\mathcal{C}}} 
\newcommand{\D}{\ensuremath{\mathcal{D}}} 
\newcommand{\E}{\ensuremath{\mathcal{E}}} 
\newcommand{\F}{\ensuremath{\mathcal{F}}} 
\newcommand{\Herb}{\ensuremath{\mathcal{H}}} 
\newcommand{\I}{\ensuremath{\mathcal{I}}} 
\newcommand{\Nat}{\ensuremath{\mathcal{N}}} 
\newcommand{\M}{\ensuremath{\mathcal{M}}} 
\newcommand{\MI}{\ensuremath{\mathcal{M}_I}} 
\newcommand{\pwst}{\ensuremath{\mathcal{P}}} 
\newcommand{\R}{\ensuremath{\mathcal{R}}} 
\newcommand{\Sch}{\ensuremath{\mathcal{S}}} 
\newcommand{\T}{\ensuremath{\mathcal{T}}} 
\newcommand{\Var}{\ensuremath{\mathcal{V}}} 
\newcommand{\Z}{\ensuremath{\mathcal{Z}}} 

\newcommand{\cart}{\ensuremath{\mbox{\textsc{cart}}}} 
\newcommand{\apl}{\ensuremath{\mbox{\textsc{apl}}}} 
\newcommand{\tit}{\ensuremath{\mbox{\textit{true}}}}
\newcommand{\tbf}{\ensuremath{\mbox{\textbf{true}}}}
\newcommand{\fit}{\ensuremath{\mbox{\textit{false}}}}
\newcommand{\fbf}{\ensuremath{\mbox{\textbf{false}}}}

\newcommand{\va}{{\tt a}}
\newcommand{\vb}{{\tt b}}
\newcommand{\vfa}{{\tt f(a)}}
\newcommand{\vfb}{{\tt f(b)}}

\newcommand{\pair}[2]{\ensuremath{\langle #1,#2 \rangle}}
\newcommand{\triple}[3]{\ensuremath{\langle #1,#2,#3 \rangle}}
\newcommand{\vc}[2]{\ensuremath{#1_0,\ldots,#1_{#2-1}}}

\newcommand{\id}{\ensuremath{\mbox{id}}}
\newcommand{\argt}{\ensuremath{\mbox{argt}}}
\newcommand{\prd}{\ensuremath{\mbox{prod}}}
\newcommand{\sq}{\ensuremath{\mbox{sq}}}
\newcommand{\dom}{\ensuremath{\mbox{dom}}}
\newcommand{\ran}{\ensuremath{\mbox{ran}}}
\newcommand{\df}{\ensuremath{\mbox{def}}}

\newcommand{\tr}{\ensuremath{\triangleright}}

\newcommand{\sol}{\ensuremath{\mbox{\textsc{sol}}}}

\newcommand{\para}{\vspace{0.05in}}

\vspace{-2mm}
\section{Introduction}
\vspace{-2mm}
\hfill
\parbox{2.0in}{
{\small
\emph{Make things as simple as possible,
but not simpler.}\\

\vspace{-.2in}
\hfill Albert~Einstein.
}
}
\vspace{3mm}

The \emph{raison d'\^{e}tre} of interval arithmetic is \emph{rigour}. 
Yet it appears that the most fundamental fact, sometimes referred to as 
the ``Fundamental Theorem of Interval Arithmetic'',
is not rigorously established.
The fact in question
can be described as follows.

Let $e$ be an expression with
$\langle x_1,\ldots,x_n \rangle$
as an ordered set of variables
(i.e. a finite sequence of distinct variables).
Let $f$ be the function in $\R^n \rightarrow \R$ that is computed by $e$.
Let the result of evaluating $e$ with intervals
$I_1,\ldots,I_n $
substituted for
$x_1,\ldots,x_n $
be an interval $Y$.
Then
\begin{equation}
\label{fact:fund}
\{f(x_1,\ldots,x_n) \mid x_1 \in I_1, \ldots, x_n \in I_n\} \subset Y
\end{equation}

Although this fact is fundamental to everything
that is done in interval arithmetic,
we have failed to find in the literature a definition
of what it means for \emph{an expression to compute a function}.
In Section~\ref{sec:remLit} we review the literature that we consulted.

In (\ref{fact:fund}) the interval $Y$ is typically considerably wider than the
range of function values. Interval analysis relies on the fact that,
as $I_1,\ldots,I_n$ become narrower, the sides in (\ref{fact:fund})
become closer to each other.
A theorem to this effect, such as 2.1.1 or 2.1.3 in \cite{nmr90}
deserves to be called Fundamental Theorem of Interval Analysis
rather than interval arithmetic.

Both theorems should rest on the foundation
provided by a definition of the function
computed by an expression.
We give such a definition for sets;
as intervals are sets,
the definition applies to intervals as a special case.

\vspace{-2mm}
\subsection{Expressions and functions}
\vspace{-2mm}
An expression is an entity consisting of \emph{symbols}\/;
it is an element of a formal language in the sense of computer science.
Some of these symbols denote operations; others are constants or variables
and denote reals or intervals,
according to the chosen interpretation.

Unlike an expression,
a numerical function is an element
of the function space $\R^n \rightarrow \R$,
for a suitable positive integer $n$.
Variables only occur in expressions,
where they can re-occur any number of times.
Variables do not occur in functions;
in fact, the notion of ``occurs in'' is not applicable to functions
in $\R^n \rightarrow \R$.
Instead, a function in
$\R^n \rightarrow \R$
is a map from $n$-tuples of reals to reals;
the elements of the $n$-tuples are properly called \emph{arguments},
rather than ``variables''.

An additional reminder of the need to distinguish between expressions and
functions is that different expressions can compute the same function.
Yet another reminder is that there exist functions that are not computable,
whereas all expressions are, like programs, instructions
for computations. Contrary to programs in general,
expressions of the type of interest to interval arithmetic
can be evaluated in finite time.
Hence the functions computed by these expressions belong to the computable
subset of functions.

Of course,
``the set of expressions'' could be made precise
by means of a formal grammar.
For the purpose of this paper,
it is sufficient to define an expression as
follows.
\vspace{-2mm}
\begin{enumerate}
\item
A variable is an expression. 
\vspace{-2mm}
\item
If $E$ is an expression and if $\varphi$ is a unary operation symbol,
then $\varphi E$ is an expression.
\vspace{-2mm}
\item
If $E_1$ and $E_2$ are expressions
and if $\diamond$ is a binary operation symbol,
then $E_1 \diamond E_2$ is an expression.
\vspace{-2mm}
\end{enumerate}
To make the definition formal,
we would have to spell out the appearance of the variables
and of the operation symbols.

In whatever way expressions are defined,
the resulting set is disjoint from the set
$\R^n \rightarrow \R$,
whatever $n$ is.
What is needed to turn
the Fundamental Fact~(\ref{fact:fund}) into a theorem
is to define ``function computed by an expression'' as a mapping
from the set of expressions in $n$ variables to
$\R^n \rightarrow \R$.
As observed above, this mapping is neither injective nor surjective.
This mapping can be called the \emph{semantics} of the language of expressions.

\vspace{-2mm}
\subsection{Remarks on the literature}
\label{sec:remLit}
\vspace{-2mm}

Moore \cite{moore66} avoids the problem of defining the function
defined by an expression by not making the distinction.
As explained in the previous section,
this is not correct.
Jaulin \emph{et. al.} \cite{jln01}, Theorem 2.2, assume that the problem
is taken care of by composition of functions,
but make unjustified simplifications.
Composition is indeed a promising approach,
which we will pursue in Section~\ref{sec:composition}.

Neumaier \cite{nmr90} does distinguish
between expressions and functions,
but the expressions as he defined them fail to be computable.
In fact, following the definition he gave in page $13$, 
every real number is an element of the set of arithmetical
expressions. 
The simplicity arises from the fact that all real numbers
are defined as (sub)expressions.
This introduces infinite expressions:
whatever notation is chosen for the reals,
most (in the sense of a subset of measure one) are infinite.
In this way effective computability is lost.

Moreover, Neumaier starts with an arithmetic expression $f$,
and then defines the interval evaluation of $f$, which he denotes
by the same symbol $f$. To deal with partial functions, 
he introduced a \emph{NaN} symbol,
and the results of operations on this symbol.
He then defined the restriction
of $f$ to its real domain $D_f=\{ x \in \R^n \mid f(x) \neq NaN\}$
to be the real evaluation of $f$.
We do not see the need for this indirect approach:
partial functions are a perfectly natural and hardly novel
generalization of functions that are total.


Ratschek and Rokne also distinguish expressions from functions.
In \cite{rtrkn88} they refer to their earlier book
\cite{rtrkn84} for a definition.
This is a mistake, because on page 23,
after a heuristic discussion of the connection between expression
and functions,
they refer to texts in logic and universal algebra for a definition.
However, these assume that all functions are total.
This is not always the case for the expressions of interest to interval
arithmetic; consider for example $\surd x$.
As only a few exotic varieties of logic allow function symbols
to be interpreted by partial functions,
it is better for interval arithmetic to use set theory
as basis for its fundamental theorems.
In fact, these exotic varieties are subject to considerable controversy
\cite{chjo91,prns93},
so not suitable as a fundament for interval analysis.

\vspace{-2mm}
\section{Set theory preliminaries}
\vspace{-2mm}
This section establishes the concepts, terminology and notation for this paper.
It is necessary because the present investigation is unusual in that all
functions are what are usually called
``partial functions''.
To avoid having to qualify with ``partial'' every time a function is mentioned,
we define ``function'' to mean what is usually referred to as 
``partial function''. In other respects, we adhere closely to standard 
expositions of set theory, such as 
\cite{hlm60,brbk39} and standard introductions such as found in
authoritative texts such as 
\cite{kll55}.

\vspace{-2mm}
\begin{definition}
A \emph{function} $f$ consists of a \emph{source}, a \emph{target}, and a \emph{map}.
The source and target are sets.
The map associates each element of a subset of the source with a unique
element of the target.
\end{definition}
\vspace{-2mm}

The set of functions with source $S$ and target $T$ is denoted by the term
$S\rightarrow T$.
If a function $f \in S \rightarrow T$ associates $x \in S$ with $y \in T$,
then one may write $y$ as $f(x)$.
When only the association under $f$ between $x$ and $y$ is relevant,
we write $x \mapsto y$.

\begin{example}
The square root is a function in $\R \rightarrow \R$
that does not associate any real
with any negative real and associates with $x \in \R$
the unique non-negative $y \in \R$ such that $y^2 = x$
if $x \geq 0$.
\end{example}

The term $f(x)$ is undefined if there is no $y \in T$
associated with $x \in S$ by $f \in S \rightarrow T$.
We take
$\{f(x) \mid x \in S\}$
to mean
$$\{y \in T \mid \exists x \in S \mbox{ and } f 
           \mbox{ associates } y \mbox{ with } x\}.$$
That is,
$\{f(x) \mid x \in S\}$
is defined even though $f(x)$ may not be defined for every
$x \in S$.

\begin{example}
$\{\surd x\mid x \in \R\}$
is defined and is the set of non-negative reals.

$\{x/y \mid x \in \{1\} \mbox{ and }  y \in \R \}$
is defined and is $\R \setminus \{0\}$.

\end{example}

The subset of $S$ consisting of $x$ with which
$f \in S \rightarrow T$ associates a $y \in T$
is called the \emph{domain} of $f$, denoted $\dom(f)$.
If $\dom(f) = S$, then $f$ is said to be a \emph{total function}.
$\{f(x) \mid x \in S\}$ is called the \emph{range} of $f$.
We introduced the unusual terms
``source'' for $S$ and ``target'' for $T$
because of the need
to distinguish them from ``domain'' and ``range''.

\vspace{-2mm}
\begin{definition}
The set of functions with source $S$ and target $T$ is denoted
$S \rightarrow T$ and is called a ``type'' or ``function space''.
\end{definition}
\vspace{-2mm}

Again, this differs from the usual meaning of $S \rightarrow T$,
where it only contains total functions.
To say that $f$ ``is of type'' 
$S \rightarrow T$ means that
$f \in S \rightarrow T$.

\vspace{-2mm}
\begin{definition}
\label{def:comp}
Let $f \in S \rightarrow T$ and $g \in T \rightarrow U$.
The \emph{composition} $g \circ f$
of $f$ and $g$ is the function in 
$S \rightarrow U$ such that
$g \circ f$ associates $x \in S$ with $z \in U$
iff there exists a $y \in T$ such that $f$
maps $x$ to $y$ and $g$ maps $y$ to $z$.
\end{definition}
\vspace{-2mm}

This is the conventional definition of composition.
It requires the target of $f$ to be the same set as the source of $g$.
Because of this requirement it is not clear
what composition Jaulin \emph{et. al.} have in mind
in \cite{jln01}, Theorem 2.2.

It follows from Definition~\ref{def:comp}
that the domain of definition of $f\circ g$
is a subset of that of $f$.

\begin{example}
$f\circ g \circ h$ has $\{0\}$ as domain if
$f \in \R \rightarrow \R$ is such that it maps $x \mapsto \surd x$,
$g \in \R \rightarrow \R$ is such that it maps $x \mapsto -x$, and
$h \in \R \rightarrow \R$ is such that it maps $x \mapsto |x|$.
In other words, $\surd (-|x|)$
is undefined for all $x \in \R$ except when $ x = 0 $.
\end{example}

Let $f \in S \rightarrow T$.
The elements of $S$ are called ``arguments'' of $f$.
Note that if a function associates an $x$ in $S$ with a $y$ in $T$,
it only so associates a \emph{single} element of $S$.
In that respect, all functions are ``single-argument'' functions.
But $S$ and $T$ may be any sets whatsoever.
Suppose $f \in \R^n \rightarrow \R$.
Now the single elements in the source of $f$, the arguments of $f$,
are $n$-tuples of reals.
Thus we interpret the usual
$f(x_1,\ldots,x_n)$
as
$f(\langle x_1,\ldots,x_n \rangle)$.

\vspace{-2mm}
\begin{definition}\label{def:cartProd}
Let
$f_1 \in S_1 \rightarrow T_1$
and
$f_2 \in S_2 \rightarrow T_2$.
The \emph{Cartesian product} of $f_1$ and $f_2$,
denoted $f_1 \times f_2$,
is a function in $S_1 \times S_2 \rightarrow T_1 \times T_2$
having domain
$
dom(f_1)\times dom(f_2)=\{\langle x_1,x_2 \rangle\mid x_1 \in \dom(f_1) \mbox{ and } x_2 \in \dom(f_2)\}
$,
and mapping
every $\langle x_1,x_2 \rangle$ in $dom(f_1)\times dom(f_2)$ to $\langle f_1(x_1),f_2(x_2) \rangle$.
\end{definition}
\vspace{-2mm}

\vspace{-2mm}
\begin{definition}
Let $f$ be a function in $S \rightarrow T$.
Let $F$ be a total function in $\pwst(S) \rightarrow \pwst(T).$
$F$ is a \emph{set extension} of $f$ iff
$\{f(x) \mid x \in X\} \subset F(X)$
for all subsets $X$ of $S$.
The total function in $\pwst(S) \rightarrow \pwst(T)$
with map $ X \mapsto \{f(x) \mid x \in X\}$
is a set extension and is called the 
\emph{canonical set extension} of $f$.
We will use $f(D)$ to denote $\{f(x) \mid x \in D\}$.
\end{definition}

\vspace{-2mm}
\Section{Intervals are sets --- interval extensions are set extensions}
\vspace{-2mm}
As we saw, partial functions have set extensions
that are total.
This is of particular interest in numerical computation,
where some important functions,
such as division and square root,
are not everywhere defined.

In some treatments of interval arithmetic
this leads to the situation
in which division by an interval containing zero is not defined.
This is not necessary:
if one regards an interval as a set and an interval extension as a set
extension,
then the interval extension is everywhere defined.
This is the approach taken in \cite{hckvnmdn01},
which will be summarized here.

A well-known fact is that the closed, connected sets of reals
have one of the following forms:
$\{x \in \R \mid x \leq b\}$,
$\{x \in \R \mid a \leq x\}$,
$\{x \in \R \mid a \leq x \leq b\}$,
as well as $\R$ itself.
Here $a$ and $b$ are reals.
Note that the empty subset of \R\ is an interval also,
as no ordering is assumed between $a$ and $b$.

The closed, connected sets of reals
are defined to be the \emph{real intervals}.
They are denoted
$[-\infty,b]$,
$[a,\infty]$,
$[a,b]$, and
$[-\infty,\infty]$.
These notations are just a shorthand for the above set expressions.
They are not meant to suggest that, for example,
$-\infty \in [-\infty,b] = \{x \in \R \mid x \leq b \}$.
This is not the case because
$[-\infty,b]$ is a set of reals and 
$-\infty$ is not a real.

The floating-point numbers are a set consisting of a finite set of reals as well
as $-\infty$ and $\infty$.
The real floating-point numbers are ordered
as among the reals. The least (greatest) element in the ordering is
$-\infty$ ($\infty$).
The floating-point intervals are the subset of the real intervals
where a bound, if it exists,
is a floating-point number.
We assume that there are at least two finite floating-point numbers.
As a result, the empty subset of \R\ is also a floating-point interval.

The floating-point intervals have the property that for every set of reals
there is a unique least floating-point interval that contains it.
This property can be expressed by means of the function $\Box$ so that
$\Box S$ is the smallest floating-point
interval containing $S \subset \R$.
Given a real number $x$, we denote by $x^-$ the greatest floating-point
number not greater than $x$, and by $x^{+}$ the least floating-point
number not less than $x$.

By themselves,
set extensions are not enough to obtain interval extensions.
They need to be used in conjunction with the function $\Box$,
as in the following definition of interval addition:
\begin{equation}
\label{def:add}
X + Y = \Box\{z \in \R \mid \exists x \in X, y \in Y\/.\/\/ x+y=z\}
\end{equation}
for all floating-point intervals $X$ and $Y$.
Compared with a definition such as 
\begin{equation}
\label{def:primAdd}
[a,b]+[c,d] = [(a+c)^-,(b+d)^+],
\end{equation}
(which is equivalent for bounded intervals)
(\ref{def:add}) has the advantage of being applicable to unbounded intervals
without having to define arithmetic operations
between real numbers and entities that are not real numbers.
Moreover, (\ref{def:add}) includes the required outward rounding.

Similarly to (\ref{def:add}) we have
\vspace{-2mm}
\begin{definition} \label{def:intvAr}
\begin{eqnarray*}
X + Y & \stackrel{\df}{=}&
     \Box\{z \in \R \mid \exists x \in X, y \in Y\/. x+y=z\}
\\
X - Y & \stackrel{\df}{=}&
     \Box\{z \in \R \mid \exists x \in X, y \in Y\/. z+y=x\}
\\
X * Y & \stackrel{\df}{=}&
     \Box\{z \in \R \mid \exists x \in X, y \in Y\/. x*y=z\}
\\
X / Y & \stackrel{\df}{=}&
     \Box\{z \in \R \mid \exists x \in X, y \in Y\/. z*y=x\}
\\
\surd X & \stackrel{\df}{=}&
     \Box\{y \in \R \mid \exists x \in X\/. y^2=x\}
\\
\end{eqnarray*}
\end{definition}

\begin{theorem}
The functions defined in Definition~\ref{def:intvAr}
map floating-point intervals to floating-point intervals,
are defined for all argument floating-point intervals,
and are set extensions of the corresponding functions from reals to reals.
\end{theorem}
This is a summary of several results in \cite{hckvnmdn01}.

\vspace{-2mm}
\begin{definition}
Let $I$ be the set of intervals.
$F \in I^n \rightarrow I$ is an \emph{interval extension}
of $f \in \R^n \rightarrow \R$ iff
$F$ is the restriction to domain $I^n \subset \R^n \rightarrow \R$
of a set extension of $f$.
$F$ is the \emph{canonical interval extension} of $f$
is defined to be $F(B) = \{f(x) \mid x \in B\}$
whenever this is an interval.
\end{definition}

\vspace{-2mm}
\Section{Semantics of expressions via set theory}
\vspace{-2mm}
\label{sec:composition}
As all but a few exotic varieties of logic restrict functions to be total,
we develop the semantics of expressions on the basis of set theory,
even though most treatments of set theory
also restrict functions to be total.
However, as we have seen,
functions in the usual set theory
are easily generalized
so that totality is not assumed.
Modifying logic so that function symbols
can be interpreted by partial functions
has graver repercussions \cite{chjo91,prns93}.

Suppose that the expression $e$ has the form $e_1+e_2$
and that $e_1$ computes $f_1 : \Rea^m \rightarrow \Rea$
and that $e_2$ computes $f_2 : \Rea^n \rightarrow \Rea$.
In such a situation, Jaulin et al. \cite{jln01} (Theorem 2.2), suggest
that the function $f$ computed by $e$ is the composition
of $+$, $f_1$ and $f_2$.

But such a composition is not possible,
as the types do not match, as required in Definition~\ref{def:comp}.
We can make a composition if we form the Cartesian product
of $f_1$ and $f_2$ and if we make additional assumptions about
$e_1$ and $e_2$.
To prepare these assumptions we need the following definition.

\vspace{-2mm}
\begin{definition}
Let $\{v_1,\ldots,v_n\}$
be the set of variables in expression $e$.
The \emph{variable sequence} of $e$
is $\langle v_1,\ldots,v_n \rangle$
if the first occurrences of the variables in $e$
are ordered according to this sequence.
\end{definition}
\vspace{-2mm}

Consider the special case where $m=n$
and where $e_1$ and $e_2$ have the same variable sequence.
Let $\delta \in \R^n \rightarrow \R^n \times \R^n$
with mapping
$$
\langle x_1,\ldots,x_n \rangle \mapsto
\langle \langle x_1,\ldots,x_n \rangle,
        \langle x_1,\ldots,x_n \rangle
\rangle
$$
As will be shown in Lemma~\ref{thm:funcExpr},
the function computed by $e_1 + e_2$ is
$+ \circ (f_1 \times f_2) \circ \delta$.
The types of $\delta$, $f_1 \times f_2$, and $+$ do match:
they are, respectively,
$\R^n \rightarrow (\R^n \times \R^n)$,
$(\R^n \times \R^n) \rightarrow \R^2$, and
$\R^2 \rightarrow \R$.
Thus it is clear the composition is defined and
that its type is $\R^n \rightarrow \R$.

But it is of course a very special case if $e_1$ and $e_2$
have the same variables in the same order of first occurrence.
To further illustrate what is needed to define a composition of
$+$, $e_1$, and $e_2$, consider another special case:
$e_1$ and $e_2$ have \emph{no} variables in common,
and their variable sequences are
$\langle v_1,\ldots,v_m \rangle$ and
$\langle w_1,\ldots,w_n \rangle$, respectively.
As will be shown in Lemma~\ref{thm:funcExpr},
the function computed by $e_1+e_2$ is again 
$+ \circ (f_1 \times f_2) \circ \delta$,
except that $\delta$ is in
$\R^{m+n} \rightarrow \R^m \times \R^n$
and has as map 
$$
\langle x_1,\ldots,x_m, y_1,\ldots,y_n  \rangle \mapsto
\langle \langle x_1,\ldots,x_m \rangle,
        \langle y_1,\ldots,y_n \rangle
\rangle
$$

Now the types of $\delta$, $f_1 \times f_2$, and $+$,
are, respectively,
$\R^{m+n} \rightarrow (\R^m \times \R^n)$,
$(\R^m \times \R^n) \rightarrow \R^2$, and
$\R^2 \rightarrow \R$.
Thus it is clear
that the composition is defined
and that its type is $\R^{m+n} \rightarrow \R$.

Finally, an example where the
subexpressions share some, but not all variables.
Consider the example where
$e_1$ is $x*y$,
$e_2$ is $y*z$,
$e$ is $e_1+e_2$,
and $\delta \in \Rea^3 \rightarrow (\Rea^2 \times \Rea^2)$
is such that $\delta$ maps
as follows:
$\triple{x_1}{x_2}{x_3} \mapsto \pair{\pair{x_1}{x_2}}{\pair{x_2}{x_3}}$
for all $x_1, x_2, x_3 \in \Rea$.
Now the functions $f_1$ and $f_2$ computed by $e_1$ and $e_2$
are the same function in $\R^2 \rightarrow \R$:
it has as map $\langle s,t \rangle \mapsto s*t$ for all reals $s$ and $t$.
Yet the function computed by $e_1+e_2$
does not have as map $s*t + s*t$:
it is a different function,
which is, however, described by the same formula
$+ \circ (f_1 \times f_2) \circ \delta$.

These three examples suggest how to define in general, for any pair
\pair{e_1}{e_2} of expressions and any domain $D$ of interpretation,
a ``distribution function'' that represents the pattern
of co-occurrences of variables in $e_1$ and $e_2$.

\vspace{-2mm}
\begin{definition}
\label{def:distrFun}
Given expressions $e_1$ and $e_2$ with variable sequences
$\langle v_1, \ldots, v_m\rangle$
and
$\langle w_1, \ldots, w_n\rangle$,
respectively.
Let $D$ be a set of values suitable for substitution
for the variables.
Let $\{i_1, \ldots, i_p\}$
and
let $\{j_1, \ldots, j_q\}$
be a partition in
$\{1,\ldots,n\}$
such that 
$\{w_{i_1}, \ldots, w_{i_p}\}$
occur in $e_1$ and
$\{w_{j_1}, \ldots, w_{j_q}\}$
do not occur in $e_1$\footnote{
Hence the variable sequence of any expression of the form
$e_1\langle \mbox{operation symbol} \rangle e_2$
is $ \langle v_1,\ldots,v_m,w_{j_1},\ldots,w_{j_q}\rangle $.

}.

The \emph{distribution function} $\delta$ for the pair
\pair{e_1}{e_2} and $D$ is the function in
$D^{m+q}\rightarrow D^m \times D^n$
that has as map
$$
\langle x_1,\ldots,x_m,y_{j_1},\ldots,y_{j_q}\rangle
\mapsto
\pair{\langle x_1,\ldots,x_m \rangle}{\langle y_1,\ldots,y_n \rangle}
$$
for all 
$x_1,\ldots,x_m, y_1,\ldots,y_n$
in $D$.

\end{definition}

\vspace{-5mm}
\begin{definition}
An \emph{interpretation} for an expression
consists of a set $D$ (the \emph{domain} of the interpretation)
and a map $M$
that maps every $n$-ary operation symbol in the expression
to a function in $D^n \rightarrow D$.

A set extension $I'$ of $I$ is said to be \emph{continuous}
if every symbol $p$ is mapped to
a continuous set extension of $M(p)$.
$I'$ is said to be \emph{canonical}
if every $n$-ary operation symbol $p$ is mapped to
a canonical set extension of $M(p)$.

\end{definition}
\vspace{-2mm}

The distribution function specifies enough of the way variables
are shared between two expressions
to support the central definition of this paper:

\vspace{-2mm}
\begin{definition}
\label{def:funcExpr}
Let $e$ be an expression and let $I$ be an interpretation
that maps each $n$-ary operation symbol in $e$ to a
function in $D^n \rightarrow D$, for $n\in \{1,2\}$.
We define by recursion on the structure of $e$,
distinguishing three cases.

Suppose $e$ is a variable.
\emph{The function computed by} $e$ \emph{under} $I$
is the identity function on $D$.

Suppose $e$ is $\varphi e_1$ where $\varphi$ is a unary operation symbol.
The function computed by $e$ \emph{under} $I$ is $f\circ f_1$,
where $f$ is the function in $D \rightarrow D$ that is the result of mapping
by $I$ of $\varphi$
and where $f_1$ is the function computed by $e_1$ under $I$.

Suppose $e$ has the form $e_1 \diamond e_2$,
where $\diamond$ is a binary operation symbol.
Suppose $\delta$ is the distribution function
for \pair{e_1}{e_2} and $D$.
The \emph{function computed by} $e$ \emph{under} $I$
is $\Diamond \circ (f_1 \times f_2) \circ \delta$,
where $\Diamond$ is the result of mapping by $I$ of $\diamond$.
\end{definition}
\vspace{-2mm}

The definition assumes that no constants occur in expressions.
We can simulate a constant by replacing it by a new variable and 
substituting the constant for that variable.
In this way the definition does not suffer a loss of generality for
expressions consisting of variables, constants, unary operators,
and binary operators.
At the expense of cumbersome notation
(or sophisticated methods to avoid this),
the function $\delta$ can be extended to cover $n$-ary operation symbols
with $n > 2$.

The definition should conform to our intuition about expression evaluation.
Suppose that $D$ is the set of integers, that the functions computed by
$e_1$ and $e_2$ yield 2 and 3, respectively.
Then the definition should ensure that the function computed by
$e_1+e_2$ yields 5 when the interpretation maps $+$ to addition over the
integers.
The following lemma confirms this intuition in general for
arbitrary binary operation symbols.

\begin{lemma}
\label{thm:funcExpr}
Let $e_1 \diamond e_2$ be the expression in Definition~\ref{def:funcExpr}.
Suppose that
$\langle a_1,\ldots,a_m\rangle \in D^m$
is substituted for
$\langle x_1,\ldots,x_m\rangle$
and that
$\langle b_1,\ldots,b_n\rangle \in D^n$
is substituted for
$\langle y_1,\ldots,y_n\rangle$.
Let 
$\langle c_1,\ldots,c_q\rangle$
be such that
$\langle a_1,\ldots,a_m,c_1,\ldots,c_q\rangle$
is substituted for
$\langle x_1,\ldots,x_m,y_1,\ldots,y_q\rangle$.

It is the case that 
\begin{eqnarray*}
& &f(\langle a_1,\ldots,a_m,c_1,\ldots,c_q\rangle)  \\
& & \;\;\;\;\; =
f_1(\langle a_1,\ldots,a_m\rangle)
\Diamond f_2(\langle b_1,\ldots,b_n\rangle),
\end{eqnarray*}
where $f$ is the function computed by $e_1 \diamond e_2$
according to Definition~\ref{def:funcExpr}.
\end{lemma}

\emph{Proof:}
\begin{eqnarray*}
f(\langle a_1,\ldots,a_m,c_1,\ldots,c_q\rangle) & = &               \\
(\Diamond \circ (f_1 \times f_2) \circ \delta)
               (\langle a_1,\ldots,a_m,c_1,\ldots,c_q\rangle) & = & \\
(\Diamond \circ (f_1 \times f_2))
        \delta(\langle a_1,\ldots,a_m,c_1,\ldots,c_q\rangle)) & = & \\
(\Diamond \circ (f_1 \times f_2))
  \langle \langle a_1,\ldots,a_m\rangle,\langle b_1,\ldots,b_n\rangle \rangle
                                                               & = & \\
\Diamond((f_1 \times f_2)
  (\langle\langle a_1,\ldots,a_m\rangle,\langle b_1,\ldots,b_n\rangle\rangle)
                                                               & = & \\
 f_1(\langle a_1,\ldots,a_m\rangle) \Diamond
 f_2(\langle b_1,\ldots,b_n\rangle).
\end{eqnarray*}

\vspace{-2mm}
\begin{lemma}
\label{ext}
Let $I$ be an interpretation for expression $e$
and let $I'$ be a set extension of $I$.
Let $f$ ($f'$) be the function computed by $e$
under the interpretation $I$ ($I'$).
Then $f'$ is a set extension of $f$.
\end{lemma}

Though a minor lemma in set theory,
the special case where the domains of $I$ and $I'$
are the reals and intervals respectively,
it plays the role of the Fundamental Theorem of Interval Arithmetic\footnote{
Except that the statement in \cite{hnsn92} inadvertently states
instead the definition of interval extension.
}.

\emph{Proof:}
We proceed by induction on the depth of the expression.
Suppose the lemma holds for all expressions of depth at most $n-1$.
Let $n$ be such that at least one of $e_1$ and $e_2$ is of depth
$n-1$ and the other is of depth at most $n-1$.
Suppose $I$ has domain $D$ and map $M$.
Let $e$ be $e_1 \diamond e_2$ and suppose that $M$ maps
$\diamond$ to 
$\Diamond$.
Let $\delta$ be the distribution function of $e_1$ and $e_2$
in that order.
Let $f_1$ and $f_2$ be the functions computed by
$e_1$ and $e_2$, respectively, under $I$.
Let $f_1'$ and $f_2'$ be the functions computed by
$e_1$ and $e_2$, respectively, under $I'$.
This gives as induction assumption that $f_1'$ and $f_2'$
are set extensions of $f_1$ and $f_2$.

Let $f$ and $f'$ be the functions computed from $e_1 \diamond e_2$
under interpretations $I$ and $I'$, respectively.
Let
$A_1,\ldots,A_m,B_1,\ldots,B_n$
be subsets of $D$ containing the elements
$a_1,\ldots,a_m,b_1,\ldots,b_n$.
Let $c_1,\ldots,c_q$ be such that $\delta$ maps
$\langle a_1,\ldots,a_m,c_1,\ldots,c_q\rangle$
to 
$\langle\langle a_1,\ldots,a_m \rangle,\langle b_1,\ldots,b_n\rangle\rangle$.

Supposing that $\Diamond '$ is a set extension of $\Diamond$,
we have
\begin{eqnarray*}
f(\langle a_1,\ldots,a_m,c_1,\ldots,c_q\rangle) & = &
\\
f_1(a_1,\ldots,a_m) \Diamond f_2(b_1,\ldots,b_n) & \in &
\\
f_1'(A_1,\ldots,A_m) \Diamond ' f_2'(B_1,\ldots,B_n) & = &
\\
f'(A_1,\ldots,A_m,C_1,\ldots,C_q), &   &
\\
\end{eqnarray*}
which is the function computed by $e$ under $I'$.
Both equalities are justified by Lemma~\ref{thm:funcExpr}.

\vspace{-2mm}
\begin{theorem}
Let $e$ be an expression with a variable sequence 
$\langle x_1, \ldots, x_n\rangle$.
Let $I$ be an interpretation for $e$, 
and $I'$ a canonical set extension of $I$.
Let $f$ ($f'$) be the function computed by $e$ 
under the interpretation $I$ ($I'$).
If each variable $x_i$ occurs only once in $e$, 
then $f'$ is the canonical set extension of $f$.
\end{theorem}

\emph{Proof:}
Following the same steps and notation as in the previous
proof, we have
$$ f'(A_1,\ldots,A_m,B_1,\ldots,B_n) =  $$
(by Lemma~\ref{thm:funcExpr})
$$ f_1'(A_1,\ldots,A_m) \Diamond ' f_2'(B_1,\ldots,B_n) = $$
(by the induction assumption)
$$ f_1(A_1,\ldots,A_m) \Diamond ' f_2(B_1,\ldots,B_n) = $$
(using that $f_1'$ and $f_2'$ are canonical set extensions
and that $e_1$ and $e_2$ have no variables in common)
\begin{eqnarray*}
\{y\in D \mid \exists y_1\in f_1(A_1,\ldots,A_m),
\\
\;\exists y_2\in f_2(B_1,\ldots,B_n). y = y_1 \Diamond y_2\} & = & 
\\
\{y\in D \mid \exists a_1\in A_1, \ldots, \exists a_m \in A_m,
\\
\;\exists b_1\in B_1, \ldots, \exists b_n \in B_n.
\\
\;y = f_1(a_1,\ldots,a_m) \Diamond f_2(b_1,\ldots,b_n) = &  &
\\
f(a_1,\ldots,a_m,b_1,\ldots,b_n)\} & = & \\ 
f(A_1,\ldots,A_m,B_1,\ldots,B_n).
\end{eqnarray*}

\vspace{-9mm}
\Section{Continuous set extensions}
\vspace{-2mm}
A fundamental fact in interval analysis
can be stated intuitively as
\begin{quote}
 We can get arbitrarily close to the range of the point evaluation
 of an expression $e$
 by computing the interval evaluation of $e$ 
 with a sufficiently narrow interval. 
 \end{quote}
So far we were only concerned with interval \emph{arithmetic}.
This fact, being a continuity property,
gets us into the realm of analysis.
So it is here that interval \emph{analysis} begins.

As the validity of the statement and proof of such a property
depends on a rigorous definition of the function computed
by an expression, it is wise to revisit the concepts and the theorems.


\vspace{-2mm}
\begin{definition}
Let $\F$ be a family of sets of $D$.
A sequence $S = \langle S_n\rangle_{n\in \Nat}$ of subsets of $D$ 
\emph{converges with respect to} $\F$ if it is nested, 
belongs to $\F$, and satisfies $\bigcap_{n\in \Nat}S_n = \{a\}$,
where $a$ is an element of $D$. 
We say that the singleton set $\{a\}$ is the \emph{limit} of $S$.
\end{definition}
\vspace{-2mm}

\begin{definition}\label{def:contSetFun}
Let $F \in \pwst(S) \rightarrow \pwst(T)$,
and let $\F_1$ and $\F_2$ be two families of sets
of $S$ and $T$, respectively.

Let $A = \langle A_n\rangle_{n\in \Nat}$
be any convergent sequence w.r.t. $\F_1$ with limit $\{a\}$.
$F$ is \emph{continuous w.r.t.} $\F_1$ and $F_2$ in $\{a\}$ iff
$\langle F(A_n)\rangle_{n\in \Nat}$
is a convergent sequence w.r.t. $\F_2$.
\end{definition}
\vspace{-2mm}

Continuity is a very strong requirement.
This raises the concern that no interesting examples might exist.
The next lemma shows that this concern is unnecessary.

\vspace{-2mm}
\begin{definition}
Let $f \in \R^{n} \rightarrow \R$, and let $\|.\|$ be the Euclidean norm
on $\R^{n}$. 
The function $f$ is \emph{Cauchy-continuous} at $c \in \R^{n}$ iff
for every $\epsilon > 0$ there exists a $\delta > 0$ such that
$\|x-c\| \leq \delta$ and $x \in dom(f)$ imply that $|f(x)-f(c)| \leq \epsilon$.

A sequence $\langle x_i \rangle_{i \in \Nat}$
with $x_i \in \R^n$ for all $i \in \Nat$
is \emph{Cauchy-convergent} to $\xi \in \R^n$ iff
for every $\epsilon > 0$ there exists an $n$ such that
$\|\xi - x_i\| \leq \epsilon$ for all $i > n$.
\end{definition}
\vspace{-2mm}

\begin{lemma}\label{lem:cauchy}
Let $f \in \R^{n} \rightarrow \R$ be Cauchy-continuous at every
$x \in \dom(f)$
and suppose $f$ has a canonical interval extension $F$.
Then $F$ is continuous w.r.t. the family of boxes of $\R^n$, 
and the family of intervals of $\R$.
\end{lemma}
\emph{Proof}\/: 
Suppose that $x$ is an element of $\R^{n}$,
and that $\langle B_n \rangle_{n\in \Nat}$
is a sequence of boxes in $\Intr^{n}$
that converges to $x$ w.r.t. the family of boxes of $\R^n$.
To prove that $F$ is continuous w.r.t. the family of boxes of $\R^n$,
and the family of intervals of $\R$, we have to show that 
the sequence $\langle F(B_n) \rangle_{n\in \Nat}$converges
w.r.t. the family of intervals of $\R$.
It is clear that this sequence is nested
and belongs to the family of intervals of $\R$.
So, we only need to show that 
$\bigcap_{n\in \Nat}F(B_n)$ is a singleton.
In fact,
$$ \bigcap_{n\in \Nat}F(B_n) = \{f(x)\}.$$
The following inclusion is obvious:
$ \{f(x)\} \subset \bigcap_{n\in \Nat}F(B_n).$
Let $y$ be an element of $\bigcap_{n\in \Nat}F(B_n)$.
This implies that for every $n \in \Nat$,
there exists $x_n$ in $B_n$ such that $f(x_n)=y$. 
Because $(B_n)_{n\in \Nat}$ is a nested sequence of boxes
that intersect in $\{x\}$, 
the sequence $(x_n)_{n\in \Nat}$ Cauchy-converges to $x$. 
Since $f$ is Cauchy-continuous at $x$, we have $f(x)=y$. 
Therefore, $ \bigcap_{n\in \Nat}F(B_n) \subset \{f(x)\},$
which proves the lemma.

\begin{lemma}\label{lem:contSetExt}
Let $f \in S \rightarrow T$, and 
let $\F_1$ and $\F_2$ be two families of sets
of $S$ and $T$, respectively. 
Let $F$ be a continuous set extension of $f$ w.r.t.
$\F_1$ and $\F_2$
and let $A = \langle A_n\rangle_{n \in \Nat}$
be a convergent sequence w.r.t. $\F_1$ with limit $\{a\}$.
Then $F(\{a\}) = \{f(a)\}$.
\end{lemma}

\emph{Proof}\/:
As $F$ is continuous w.r.t.  $\F_1$ and $\F_2$, 
$\langle F(A_n)\rangle_{n\in \Nat}$
is a convergent sequence w.r.t. $\F_2$ with limit, say, $\{b\}$.
As $F$ is a set extension of $f$ we have that
$\{f(x)\mid x \in A_i\} \subset F(A_i)$, for all $i \in \Nat$.
As $a \in A_i$ for all $i \in \Nat$,
we have that
$f(a) \in \{f(x)\mid x \in A_i\}$ for all $i \in \Nat$.
Hence $f(a) \in \bigcap_{i \in \Nat} F(A_i) = \{b\}$.
So we must have $f(a) = b$.

We are interested in interval extensions that are not canonical,
yet are continuous.

Starting from a family $\F$ of sets of a set $D$, 
we can construct a family of sets $\F_n$ of $D^n$, 
for any natural number $n$, 
by taking all the Cartesian products of any $n$ sets in $\F$.
So, for any natural number $n$, 
and for any function $F \in \pwst(D^n) \rightarrow \pwst(D)$, 
we can study the continuity of $F$ w.r.t. $\F_n$ 
that was constructed from $\F$. 
In this way, we treat the continuity of $F$ 
by referring to $\F$ instead of $\F_n$.

In what follows, we suppose that
the family of sets $\F$ of the domain $D$ of an interpretation
is given, and that the continuity of a set extension of
an $n$-ary operation is based on this family. 
So, we will not use ``w.r.t.'' from now on.
In the case where $D$ is $\R$, 
$\F$ is the family of intervals in $\R$.

\vspace{-2mm}
\begin{definition}
Let $I$ be an interpretation with domain $D$ and map $M$.
A set extension $I'$ of $I$ is said to be \emph{continuous}
if every symbol $p$ is mapped to
a continuous set extension of $M(p)$.
$I'$ is said to be a \emph{canonical interval extension} of
$I$ iff every symbol $p$ is mapped to a canonical interval
extension of $M(p)$.
\end{definition}
\vspace{-2mm}

\begin{theorem}\label{theo:contSetExt}
Let $e$ be an expression. 
Let $I$ be an interpretation for $e$, 
and let $I'$ be a continuous set extension of $I$.
Let $f$ ($F$) be the function computed by $e$ 
under the interpretation $I$ ($I'$).
Then $F$ is a continuous set extension of $f$.
\end{theorem}

\emph{Proof:}
From Lemma~\ref{ext}, the function $F$ is a set extension of $f$.
So we only need to prove that $F$ is continuous.
To do so,
we proceed by induction on the depth of the expression $e$.
The theorem holds when $e$ has no subexpressions,
that is, when $e$ is a variable.
In that case $f$ and $F$ are the identity functions,
independently of $I$ and $I'$.
The identity function in $\pwst(D) \rightarrow \pwst(D)$
is continuous.

This takes care of the base of the inductive proof.
Let the induction assumption be
that the theorem holds for all expressions of depth at most $d-1$.
Let $e$ be the expression $e_1 \diamond e_2$, 
where one of the subexpressions
has depth $d-1$ and the other has depth at most $d-1$.
Suppose that the interpretation $I$ has domain $D$
and maps $\diamond$ to $\Diamond$.
Let the interpretation $I'$ have $\pwst(D)$ as domain
and let it map $\diamond$ to $\Diamond '$,
a continuous set extension of $\Diamond$. 
Let $\delta$ be the distribution function with $D$ for $e_1$ and $e_2$
in that order.
Let $c_1,\ldots,c_q$ be such that $\delta$ maps
$\langle a_1,\ldots,a_m,c_1,\ldots,c_q\rangle$
to 
$\langle\langle a_1,\ldots,a_m \rangle,\langle b_1,\ldots,b_n\rangle\rangle$.

Let $F$, $F_1$, and $F_2$ be the functions
computed under $I'$ by $e$, $e_1$, and $e_2$, respectively.
Suppose that 
$\langle A^i_1 \rangle_{i \in \Nat},\ldots,\langle A^i_m \rangle_{i \in \Nat}$
and
$\langle B^i_1 \rangle_{i \in \Nat},\ldots,B^i_n \rangle_{i \in \Nat}$
are sequences of subsets of $D$ that converge respectively to 
$\{a_1\},\ldots,\{a_m\}$ and
$\{b_1\},\ldots,\{b_n\}$.
According to the induction assumption
$F_1$ and $F_2$ are continuous set extensions.
This implies that
$\langle \langle F_1(A^i_1,\ldots,A^i_m),
            F_2(B^i_1,\ldots,B^i_n)\rangle \rangle_{i \in \Nat}$
converges to $\{\langle f_1(a_1,\ldots,a_m), f_2(b_1,\ldots,b_n)\rangle\}$,
by Lemma~\ref{lem:contSetExt}.

Let
$\langle C^i\rangle_{i \in \Nat}$ be any such that
$\delta(C^i) =
\langle \langle A^i_1, \ldots, A^i_m \rangle,
        \langle B^i_1, \ldots, B^i_n \rangle
\rangle
$
and such that $\langle C^i\rangle_{i \in \Nat}$ converges to 
$\{
\langle \langle a_1, \ldots, a_m \rangle,
        \langle b_1, \ldots, b_n \rangle
\}$

We show that $F$ is a continuous set extension of $f$ by showing
that $\langle F(C^i)\rangle_{i \in \Nat}$
converges to $\{f(\langle a_1, \ldots, a_m,c_1, \ldots, c_q \rangle)\}$.
To do so, we need to show that the sequence $\langle F(C^i)\rangle_{i \in \Nat}$
is nested, and that $\bigcap_{i \in \Nat} F(C_i)$ is
the right value, namely $\{f(\langle a_1, \ldots, a_m,c_1, \ldots, c_q \rangle)\}$.
$$
F(C^{i+1}) =
$$
(by Definition~\ref{def:funcExpr})
$$
(\Diamond ' \circ (F_1 \times F_2)\circ \delta)(C^{i+1}) =
$$
(by application of $\delta$)
$$
(\Diamond ' \circ (F_1 \times F_2))
          (\langle \langle A^{i+1}_1, \ldots, A^{i+1}_m \rangle,
                   \langle B^{i+1}_1, \ldots, B^{i+1}_n \rangle \rangle) =
$$
(by Definition~\ref{def:cartProd})
$$
\Diamond '
( \langle
        F_1(\langle A^{i+1}_1, \ldots, A^{i+1}_m \rangle),
        F_2(\langle B^{i+1}_1, \ldots, B^{i+1}_m \rangle)
\rangle) \subset
$$
(by the induction assumption and continuity of $\Diamond '$)
$$
\Diamond '
( \langle
        F_1(\langle A^{i}_1, \ldots, A^{i}_m \rangle),
        F_2(\langle B^{i}_1, \ldots, B^{i}_m \rangle)
\rangle) =F(C^i),
$$
which proves that $\langle F(C^i)\rangle_{i \in \Nat}$ is nested.
As for the convergence to the right value, we observe the following:

$$
\bigcap_{i \in \Nat} F(C^i) =
$$
(by Definition~\ref{def:funcExpr})
$$
\bigcap_{i \in \Nat} (\Diamond ' \circ (F_1 \times F_2)\circ \delta)(C^i) =
$$
(by application of $\delta$)
$$
\bigcap_{i \in \Nat} (\Diamond ' \circ (F_1 \times F_2))
          (\langle \langle A^i_1, \ldots, A^i_m \rangle,
                 \langle B^i_1, \ldots, B^i_n \rangle \rangle) =
$$
(by Definition~\ref{def:cartProd})
$$
\bigcap_{i \in \Nat} \Diamond '
( \langle
        F_1(\langle A^i_1, \ldots, A^i_m \rangle),
        F_2(\langle B^i_1, \ldots, B^i_m \rangle)
\rangle) =
$$
(by continuity of $\Diamond '$)
$$
\Diamond '
( \langle
  \bigcap_{i \in \Nat}F_1(\langle A^i_1, \ldots, A^i_m \rangle), \\
  \bigcap_{i \in \Nat}F_2(\langle B^i_1, \ldots, B^i_m \rangle)
\rangle) =
$$
(by the induction assumption)
$$
\Diamond '
( \langle
  \{f_1(\langle a_1, \ldots, a_m \rangle)\},
  \{f_2(\langle b_1, \ldots, b_n \rangle)\}
\rangle) =
$$
(by Lemma~\ref{lem:contSetExt})
$$
  \{f_1(\langle a_1, \ldots, a_m \rangle)\} \Diamond
  \{f_2(\langle b_1, \ldots, b_n \rangle)\} =
$$
(because $f$ is the function computed by $e_1\diamond e_2$)
$$
  \{f(\langle a_1, \ldots, a_m,c_1, \ldots, c_q \rangle)\},
$$
which shows that $F = \Diamond '\circ(F_1 \times F_2)\circ \delta$
is a continuous set extension of $f$, the function computed by $e$.

\begin{corollary}
Let $f \in \R^n \rightarrow \R$ be the function
computed by an expression $e$ under an interpretation $I$
that assigns Cauchy-continuous functions
to the operation symbols in $e$.
Let $F$ be the function computed by $e$ under the canonical
interval extension of $I$.
Let $\langle A_i\rangle_{i\in \Nat}$
be nested boxes converging to $\{a\}$.
Then $\langle F(A_i)\rangle_{i\in \Nat}$ 
is a sequence of nested intervals converging to $\{f(a)\}$.
\end{corollary}
In interval analysis, this corollary plays the role of Fundamental Theorem.

\emph{Proof:}
Since the image of any box by a Cauchy-continuous function
is an interval,
the interval extension associated with each operation
symbol is canonical
(every Cauchy-continuous function has a canonical interval extension).
Using Lemma~\ref{lem:cauchy}, these interval extensions are continuous.
By Theorem~\ref{theo:contSetExt}, $F$ is continuous.
By Definition~\ref{def:contSetFun},
$\langle F(A_i)\rangle_{i\in \Nat}$ converges to $\{f(a)\}$.

\vspace{-2mm}
\section{Conclusions}
\vspace{-2mm}


The fact that the result
of an expression evaluation in intervals gives a result
that contains the range of values of the function
computed by the expression
cannot be a mathematical theorem
without a mathematical definition
of what it means for a function to be computed by an expression.
In this paper we give such a definition
and prove the theorem on the basis of it.

Another fundamental assumption
in the use of intervals is that,
as we make the intervals
in an interval evaluation of an expression narrower,
the interval result gets closer
to the range of values of the function
computed by the expression.
We use our definition to prove a theorem to this effect.

Our starting point in all this is that intervals are sets
and that, therefore,
interval extensions of functions are set extensions of functions.
The latter concept is an old one in set theory
and is more widely applicable.
Our definition and two main theorems are stated in terms of sets,
so apply to intervals as special cases.

This is of course only of interest to those who believe in sets
as foundation of mathematics. A radically different approach to the
fundamental theorems of interval analysis is found in Paul Taylor's
work (see for example \cite{TaylorP:intawi}).
Here the starting point is topology, axiomatically founded rather
than set-theoretically.

If it seems that our proposed foundations
for interval methods are overly complex
in comparison with the way they are given in the literature,
we are comforted by Einstein's dictum:
\emph{Make things as simple as possible, but not simpler.}

\vspace{-2mm}
\section{Acknowledgements}
\vspace{-5mm}
This research was supported by the University of Victoria and by
the Natural Science and Engineering Research Council of Canada.
We owe a great debt of gratitude to our anonymous reviewer whose
extremely detailed and helpful report has helped us to improve
this paper.
\vspace{-2mm}

\end{document}